\documentclass{amsart}

\usepackage{amsmath,amsfonts,amssymb,amsthm,latexsym}
\usepackage{mathtools}

\usepackage[T1]{fontenc}
\usepackage[USenglish]{babel}
\usepackage{graphicx} 
\usepackage{epstopdf}
\usepackage{tikz}
\usetikzlibrary{matrix,arrows,cd}

\usepackage[all]{xy}
\usepackage{hyperref}

\usepackage{url}

\newtheorem{thm}{Theorem}[section]
\newtheorem{prop}[thm]{Proposition}
\newtheorem{lem}[thm]{Lemma}

\theoremstyle{definition}
\newtheorem{defi}[thm]{Definition} 
\newtheorem{example}[thm]{Example} 

\newtheorem{remark}[thm]{Remark}

\numberwithin{equation}{section}

\DeclareMathOperator{\Id}{Id}
\DeclareMathOperator{\car}{char}
 \DeclareMathOperator{\Ad}{Ad}
 \DeclareMathOperator{\Conj}{Cj}
  \DeclareMathOperator{\Ob}{Ob}

\begin{document}

	\title[Braiding I]{Braiding for categorical algebras and crossed modules of algebras I: \\
		Associative and Lie algebras}
	\author{A. Fern\'andez-Fari\~na}
\address{[A. Fern\'andez-Fari\~na] Department of Matem\'aticas, University of Santiago de Compostela, 15782, Spain.}
\email{alejandrofernandez.farina@usc.es}

	\author{M. Ladra}
\address{[M. Ladra] Department of Matem\'aticas, Institute of Matem\'aticas, University of Santiago de Compostela, 15782, Spain.}
\email{manuel.ladra@usc.es}

\thanks{This work was partially supported by Agencia Estatal de Investigaci\'on (Spain), grant MTM2016-79661-P (European FEDER
support included, UE). The first author is also supported by a scholarship of Xunta de Galicia (Spain), grant ED481A-2017/064, Xunta de Galicia (Spain).}

\begin{abstract}
In this paper we study the categories of braided categorical associative algebras and braided crossed modules of associative algebras
and we relate these structures with the categories of braided categorical Lie algebras and braided crossed modules of Lie algebras.

\end{abstract}

\maketitle
\section*{Introduction}\addcontentsline{toc}{section}{Introduction}

Monoidal categories were introduced by Jean B\'enabou~\cite{Ben63} and  Saunders Mac Lane~\cite{MLane63} in order to generalize the idea of tensor product in arbitrary categories.

It is well know that, in the case of the usual tensor product for vector spaces, there is a natural isomorphism between $V\otimes W$ and $W\otimes V$.
 In order to study if this property also holds in an arbitrary monoidal category, i.e.\ when the tensor product is commutative, Joyal and Street defined in~\cite{JSBM}
  the concept of braiding for monoidal categories as a natural isomorphism $\tau_{A,B}\colon A\otimes B\xrightarrow{}B\otimes A$.

When we try to study the concept of braiding for the simplest case of monoidal categories (categorical monoids or internal categories in the category of monoids)
we encounter the obstacle that not all internal morphisms are internal isomorphisms, so the braiding cannot be an arbitrary internal morphism verifying simple properties.
To avoid this problem we can work with (strict) categorical groups instead of categorical monoids, obtaining immediately the definition of braided categorical group (see~\cite{JSBM}).

On the other hand, in 1949 Whitehead~\cite{Whi} introduced the notion of crossed module of groups as an algebraic model for 2-type homotopy spaces
(i.e.\ connected spaces with trivial  homotopy groups in dimension $>2$). In 1984, Conduch\'e~\cite{Condu} introduced the notion of braided crossed module of groups
 as a particular case of $2$-crossed module of groups.

It is well know that the categories of crossed modules of groups and categorical groups are equivalent, and Joyal and Street proved in~\cite{JSBM}
 that the notion of braiding for categorical groups gives an equivalent category to the category of braided crossed modules of groups introduced by Conduch\'e~\cite{Condu}.

The notions of crossed modules for associative algebras~\cite{DeLu66}, Lie algebras~\cite{KassLod} and Leibniz algebras~\cite{LodayPira} appear
trying to imitate crossed modules of groups, and it was proven that the correspondent categories are equivalent to their respective internal categories.

Keeping in mind what is done for groups, in this paper we will give definitions of braidings for the aforementioned internal categories and crossed modules.
The case of associative algebras is not complex, because the associativity allows us to work in a natural way with braidings on semigroupal categories~\cite{CrYet98}.
 The notion of braiding for Lie algebras was already given by Ulualan~\cite{Ulua}.
On the other hand, Ellis~\cite{Ellis93} defined the notion of $2$-crossed module of Lie algebras, also studied by Martins and Picken~\cite{M&P}.
We will use a slightly different definition for braiding for crossed modules of Lie algebras than the one given by Ulualan~\cite{Ulua},
 since we want a parallelism between the examples of braided crossed modules of groups and braided crossed modules of Lie algebras and
 we want braided crossed modules to be a particular case of $2$-crossed modules, as in happens in groups case.
The Leibniz algebras case will be studied in the second part of this paper~\cite{FFBraidII}.

This first manuscript is organized as follows.
In the preliminaries we will recall some basic definitions and we give the notion of braiding for semigroupal categories.
In Section~\ref{S:braidgrp} we will show that, in all the internal categories we will work with, all internal morphisms are internal isomorphisms;
which motivates the study of braided crossed modules of groups instead of braided categorical monoids.
In Section~\ref{S:braidassalg} we introduce the notions of braided categorical associative algebra and braided crossed module of associative algebras
 and we show that these categories, as in the groups case, are equivalent.
In Section~\ref{S:braidLiealg} we motivate the definition given by Ulualan~\cite{Ulua} for braided crossed modules of Lie algebras using our definition
of braiding for crossed modules of associative algebras and we give a simpler definition when $\car(K)\neq 2$.
We will also discuss about a different definition of braided crossed module of Lie algebras showing its relationship with the associative case.
In Section~\ref{S:Lienonabtensor} we see the non-abelian tensor product of groups as an example of braided crossed module of groups.
 Moreover, with our definition of braiding for crossed modules of Lie algebras, we obtain in a parallel way an example of braiding using the non-abelian tensor product.

\section{Preliminaries}

\subsection{Internal Categories}

\begin{defi}
	Let \textit{\textbf{C}} be a category with pullbacks.
	
	An \emph{internal category} in $\textit{\textbf{C}}$ consist of two objects $C_1$ (\emph{morphisms object}) and $C_0$ (\emph{objects object}) of $\textit{\textbf{C}}$, together with the four following morphisms:
	\[
	\xymatrix@=3em{  C_0  \ar[r]|-{e}     & C_1\ar@<1ex>[l]^-{s} \ar@<-1ex>[l]_-{t}  & C_1\times_{C_0}C_1 \ar[l]_-{k}, }
	\]
	where $C_1\times_{C_0}C_1$ is the pullback of $t$ and $s$.
	
	 $s$ is called \emph{source morphism}, $t$ is called \emph{target morphism}, $e$ is called \emph{identity mapping morphism} and $k$ is called \emph{composition morphism}.

In addition, the morphisms must verify commutative diagrams that express the usual category laws (see~\cite{Baez04}).

	If the conditions are allowed we will refer to the internal category by the $6$-tuple $(C_1,C_0,s,t,e,k)$.
\end{defi}

\begin{defi}
	Let $\mathcal{C}=(C_1,C_0,s,t,e,k)$ and $\mathcal{C'}=(C_1',C_0',s',t',e',k')$ be two internal categories in \textbf{\textit{C}}.
	
	An \emph{internal functor} is a pair of morphisms $(F_1,F_0)$, with $F_1\colon C_1\xrightarrow{} C_1'$ and $F_0\colon C_0\xrightarrow{} C_0'$
 such that must verify commutative diagrams corresponding to the usual laws satisfied by a functor (see~\cite{Baez04}).

	We will denote the internal functor by $(F_1,F_0)\colon\mathcal{C}\xrightarrow{}\mathcal{C'}$.
\end{defi}

Composition of internal functors is defined in the obvious way. This allows us to construct the category of internal categories and internal functors
 in a category with pullbacks \textbf{\textit{C}}, denoted by $\textbf{\textit{ICat}}(\textbf{\textit{C}})$.

An internal category in $\textbf{\textit{C}}$ will be also called categorical object in $\textbf{\textit{C}}$.

\subsection{Algebras}

\begin{defi}
	Let $K$ be a field and $(M,*)$ be a $K$-algebra.
	
	A \emph{derivation} over $(M,*)$ is a $K$-linear map $D\colon M\xrightarrow{}M$ verifying the \emph{Leibniz rule}:
	\begin{align*}
	D(x*y)=D(x)*y+x*D(y), \ x,y\in M.
	\end{align*}
\end{defi}

\begin{remark}
	Let $(M,*)$ be a $K$-algebra.
	It is immediate to check that, if we take $x\in M$, the map $R(x)\colon M\xrightarrow{}M$ defined by $R(x)(y)=y*x$ (right multiplication) is $K$-linear using the $K$-bilinearity.
\end{remark}

\begin{defi}
	We will say that the $K$-algebra $(M,*)$ is a \emph{(right) Leibniz $K$-algebra} if and only if $R(x)$ is a derivation over $(M,*)$ for all $x\in M$.
We denote $x*y=:[x,y]$ and call the operation $[-,-]$ \emph{Leibniz bracket}.
	
	If, in addition, $(M,[-,-])$ is an alternate $K$-algebra ($[x,x]=0, \  x\in M$) we will say that it is a \emph{Lie $K$-algebra} and we will call the operation $[-,-]$ \emph{Lie bracket}.
\end{defi}

\begin{remark}
	 The fact that $R(z)$ is a derivation for all $z\in M$ can be seen in the next identity for $x,y,z\in M$, called the \emph{Leibniz identity}:
	\begin{equation*}
	[x,[y,z]]=[[x,y],z]-[[x,z],y].
	\end{equation*}
	If in addition the $K$-algebra is anticommutative (for example the Lie $K$-algebras), we can rewrite the equality, obtaining the \emph{Jacobi identity}:
	\begin{equation*}
	[x,[y,z]]+[y,[z,x]]+[z,[x,y]]=0.
	\end{equation*}
\end{remark}

\begin{defi}
	If $(M,*)$ and $(N,\star)$ are $K$-algebras, an \emph{homomorphism} between them is a $K$-linear map $M\xrightarrow{f}N$ such that $f(x*y)=f(x)\star f(y)$.
\end{defi}

	We have the categories $\textbf{\textit{AssAlg}}_K$, $\textbf{\textit{LieAlg}}_K$ and $\textbf{\textit{LeibAlg}}_K$ taking, as objects (respectively),
 associative, Lie or Leibniz $K$-algebras, and as morphisms, homomorphisms of $K$-algebras between them.	
	
	We will denote by $\textbf{\textit{Vect}}_K$ and $\textbf{\textit{Grp}}$ the categories of $K$-vector spaces and groups.
	
\subsection{Crossed Modules}
\subsubsection{Crossed Modules of groups}\hfill

The crossed modules of groups were introduced by Whitehead in~\cite{Whi}.

\begin{defi}
	A \emph{crossed module of groups} is a $4$-tuple $(G,H,\cdot,\partial)$ where $G$ and $H$ are groups, $\cdot$ is an action of $H$ on $G$ by automorphisms,
 $\partial\colon G\xrightarrow{}H$ is a group homomorphism and the following properties are satisfied for $g,g'\in G$, $h\in H$:
	
	$\partial$ is $H$-equivariant map (we suppose the conjugation action of $H$ on itself), i.e.
	\begin{equation}\label{XGr1}
	\partial(h\cdot g)=h\partial(g)h^{-1}.\tag{XGr1}
	\end{equation}
	
	\emph{Peiffer identity}:
	\begin{equation}\label{XGr2}
	\partial(g)\cdot g'=gg'g^{-1}.\tag{XGr2}
	\end{equation}
\end{defi}

\begin{example}
	It is clear from the definitions that if $G$ is a group, then $(G,G,\Conj,\Id_G)$ is a crossed module of groups, where $\Conj$ is the conjugation action.
\end{example}

\begin{defi}
	An \emph{homomorphism of crossed modules of groups} between $(G,H,\cdot,\partial)$ and $(G',H',*,\partial')$
	 is a pair of group homomorphisms, $f_1\colon G \xrightarrow{}G'$ and $f_2\colon H\xrightarrow{}H'$ such that:
	\begin{align*}
	f_1(h\cdot g)&=f_2(h)*f_1(g), \quad  g\in G, h\in H, \tag{XGrH1}\\
	\partial' \circ f_1&=f_2\circ \partial.\tag{XGrH2}
	\end{align*}	
\end{defi}

\begin{remark}
	There is and equivalence between the categories $\textbf{\textit{ICat}}(\textbf{\textit{Grp}})$ and crossed modules of groups (see~\cite{Baez04}).
\end{remark}

\subsubsection{Crossed Modules of associative algebras}\hfill

The definition of action in the associative algebras case is the following one.

\begin{defi}
	Let $N$ and $M$ be two associative $K$-algebras.
	
	An \emph{associative action} of $N$ on $M$ is a pair of $K$-bilinear maps $*=(*_1,*_2)$, where $*_1\colon N\times M\xrightarrow{}M$ and $*_2\colon M\times N\xrightarrow{}M$, $(n,m)\mapsto n*_1 m$ and $(m,n)\mapsto m*_2 n$, verifying for $n\in N$ and $m\in M$:
	\begin{align}
	n*_1 (mm')&=(n*_1 m)m',\tag{AAs1}\\
	n*_1 (m*_2 n')&=(n*_1 m)*_2 n',\tag{AAs2}\\
	n*_1(n'*_1m)&=(nn')*_1 m,\tag{AAs3}\\
	m*_2(nn')&=(m*_2 n)*_2 n',\tag{AAs4}\\
	m(n*_1m')&=(m*_2 n)m',\tag{AAs5}\\
	m(m'*_2n)&=(mm')*_2 n.\tag{AAs6}
	\end{align}
	for $m,m'\in M$, $n,n'\in N$.
\end{defi}

\begin{remark}
	If we change the notation of $*_1$ and $*_2$ by $*$ in both  cases, where $*$ is the multiplication,
 the axioms of the associative actions are all possible rewrites of the associativity when we choose two elements in $M$ and one in $N$ or one in $M$ and two $N$.
	
	In particular, if $*$ is the multiplication of an associative $K$-algebra, we have that the pair $(*,*)$ is an associative action of $M$ on itself.
\end{remark}

The definition of crossed module of associative algebras was given by Dedecker and Lue in~\cite{DeLu66}.

\begin{defi}
	A \emph{crossed module of associative $K$-algebras} is a $4$-tuple $(M,N,*,\partial)$ where $M$ and $N$ are associative $K$-algebras, $*=(*_1,*_2)$
 is an associative action of $N$ on $M$, $\partial\colon M\xrightarrow{}N$ is an associative $K$-homomorphism and the following properties are verified for $m,m'\in M$, $n\in N$:
	
	$\partial$ is an $N$-equivariant associative $K$-homomorphism (we suppose the action of $N$ on itself is the product), i.e.
	\begin{equation}\label{XAs1}
	\partial(n*_1 m)=n\partial(m) \quad  \text{and} \quad \partial(m*_2 n)=\partial(m)n.\tag{XAs1}
	\end{equation}
	
	Peiffer identity:
	\begin{equation}\label{XAs2}
	\partial(m)*_1 m'=mm'=m*_2\partial(m').\tag{XAs2}
	\end{equation}
\end{defi}

\begin{example}
	If $M$ is an associative $K$-algebra then $(M,M,(*,*),\Id_M)$ is a crossed module of associative $K$-algebras.
\end{example}

\begin{defi}
	An \emph{homomorphism of crossed modules of associative $K$-algebras} between $(M,N,\cdot,\partial)$ and $(M',N',*,\partial')$
	is a pair of associative $K$-homomorphisms, $f_1\colon M \xrightarrow{}M'$ and $f_2\colon N\xrightarrow{}N'$ such that:
	\begin{align*}
	f_1(n\cdot_1 m)=f_2(n)*_1f_1(m) \quad  &\text{and} \quad  f_1(m\cdot_2 n)=f_1(m)*_2 f_2(n),\tag{XAssH1}\\
	\partial' \circ f_1&=f_2\circ \partial,\tag{XAssH2}
	\end{align*}
	for $m\in M$, $n\in N$.	
\end{defi}

	We will denote by $\textbf{\textit{X}}(\textbf{\textit{AssAlg}}_K)$ the category of crossed modules of associative $K$-algebras and its homomorphisms.

\begin{remark}
	As in the case of groups we have an equivalence between the categories $\textbf{\textit{ICat}}(\textbf{\textit{AssAlg}}_K)$ and $\textbf{\textit{X}}(\textbf{\textit{AssAlg}}_K)$.
A proof can be found in~\cite{ThRa}.
\end{remark}

\subsubsection{Crossed Modules of Lie algebras}\hfill

We have an analogous definition for the case of Lie $K$-algebras.
Crossed modules of Lie $K$-algebras was introduced by Kassel and Loday in~\cite{KassLod}.

\begin{defi}
	Let $M$ and $N$ two Lie $K$-algebras, a \emph{Lie (left-)action of $N$ on $M$} is a $K$-bilinear map $\cdot\colon N\times M\longrightarrow M$, $(n,m)\longmapsto n\cdot m$, satisfying:
	\begin{align}
	[n,n']\cdot m&=n\cdot(n'\cdot m)-n'\cdot(n\cdot m),\tag{ALie1}\label{ALie1}\\
	n\cdot[m,m']&=[n\cdot m,m']+[m,n\cdot m'],\tag{ALie2}\label{ALie2}
	\end{align}
	for $n,n'\in N$, $m,m'\in M$.
\end{defi}

\begin{example}
	Note that the two identities are, if we denote $\cdot=[-,-]$, the two possible rewrites of the Jacobi identity taking two elements in $N$ or two in $M$.
	
	In particular, we have that the adjoint map $\Ad(x)\colon M\xrightarrow{}M$, where $M$ is a Lie $K$-algebra and $x\in M$, defined by $\Ad(x)(y)=[x,y]$, is a Lie action of $M$ on itself.
\end{example}

\begin{defi}
	A \emph{crossed module of Lie $K$-algebras} is a $4$-tuple $(M,N,\cdot,\partial)$ where $M$ and $N$ are Lie $K$-algebras, $\cdot$ is a Lie action of $N$ on $M$,
 $M\xrightarrow{\partial}N$ is a Lie $K$-homomorphism and the following properties are satisfied for $n\in N$ and $m,m'\in M$:
	
	$\partial$ is $N$-equivariant Lie $K$-homomorphism (we suppose the adjoint action of $N$ on itself), i.e.
	\begin{equation}\label{XLie1}
	\partial(n\cdot m)=[n,\partial(m)].\tag{XLie1}
	\end{equation}
	
	Peiffer identity:
	\begin{equation}\label{XLie2}
	\partial(m)\cdot m'=[m,m'].\tag{XLie2}
	\end{equation}
\end{defi}

\begin{example}
	As in the case of groups we have the example of crossed module of Lie $K$-algebras $(M,M,[-,-],\Id_M)$ where $M$ is a Lie $K$-algebra. Note that the adjoint action replaces the conjugation action in this case.
\end{example}

\begin{defi}
	An \emph{homomorphism of crossed modules of Lie $K$-algebras} between $(M,N,\cdot,\partial)$ and $(M',N',*,\partial')$
	is a pair of Lie $K$-homomorphisms, $f_1\colon M \xrightarrow{}M'$ and $f_2\colon N\xrightarrow{}N'$ such that:
	\begin{align*}
	f_1(n\cdot m)&=f_2(n)*f_1(m),\quad  m\in M, n\in N,\tag{XLieH1}\\
	\partial' \circ f_1&=f_2\circ \partial.\tag{XLieH2}
	\end{align*}	
\end{defi}

We want to show that there is a natural way to connect the crossed modules of associative $K$-algebras and the crossed modules of Lie $K$-algebras.

\begin{prop}\label{AsAc->LieAc}
	Let $M$ and $N$ be two associative $K$-algebras.
	
	We will denote the Lie $K$-algebra associated to an associative $K$-algebra $A$ as $A^\mathcal{L}$. That means $A^\mathcal{L}$ is the Lie $K$-algebra with the operation $[a,a']=aa'-a'a$.
	
	Then, if $*=(*_1,*_2)$ is an associative action of $N$ on $M$, we have that the map $[-,-]_*\colon N\times M\xrightarrow{}M$,
defined as $(n,m)\mapsto [n,m]_*=n*_1 m-m*_2 n$, is a Lie action between $N^\mathcal{L}$ and $M^\mathcal{L}$.
	\begin{proof}
		Let $n,n'\in N^\mathcal{L}$, $m,m'\in M^\mathcal{L}$.
		In first place we will prove~\eqref{ALie1}.
		\begin{align*}
		&[n,[n',m]_*]_*-[n',[n,m]_*]_*\\
		&=n*_1(n'*_1 m)-n*_1(m*_2 n')-(n'*_1 m)*_2 n+(m*_2 n')*_2 n\\
		&-n'*_1(n*_1 m)+n'*_1(m*_2 n)+(n*_1 m)*_2 n'-(m*_2 n)*_2 n'\\
		&=(nn')*_1 m-n*_1(m*_2 n')-n'*_1 (m*_2 n)+m*_2 (n'n)\\
		&-(n'n)*_1 m+n'*_1(m*_2 n)+n*_1 (m*_2 n')-m*_2 (nn')\\
		&=(nn')*_1 m+m*_2 (n'n)-(n'n)*_1 m-m*_2 (nn')\\
		&=[n,n']*_1 m-m*_2 [n,n']=[[n,n'],m]_*.
		\end{align*}
		We have~\eqref{ALie2} too:
		\begin{align*}
		&[[n,m]_*,m']+[m,[n,m']_*]\\
		&=(n*_1 m)m'-m'(n*_1 m)-(m*_2 n)m'+m'(m*_2 n)\\
		&+m(n*_1m')-(n*_1m')m-m(m'*_2n)+(m'*_2n)m\\
		&=n*_1 (mm')-(m'*_2n)m-(m*_2 n)m'+(m'm)*_2 n\\
		&+(m*_2n)m'-n*_1(m'm)-(mm')*_2n+(m'*_2n)m\\
		&=n*_1 (mm')+(m'm)*_2 n-n*_1(m'm)-(mm')*_2n\\
		&n*_1[m,m']-[m,m']*_2 n=[n,[m,m']]_*.
		\end{align*}
		In both we use the ``associativity'' of associative actions.
	\end{proof}
\end{prop}

\begin{prop}
	If $(M,N,*,\partial)$ is a crossed module of associative $K$-algebras, then  $(M^\mathcal{L},N^\mathcal{L},[-,-]_*,\partial)$ is a crossed module of Lie $K$-algebras.
	\begin{proof}
		Is immediate that $\partial$ is a Lie $K$-homomorphism.
		
		Let $m,m'\in M^\mathcal{L}$, $n\in N^\mathcal{L}$.
		\begin{align*}
		&\partial([n,m]_*)=\partial(n*_1 m)-\partial(m*_2 n)=n\partial(m)-\partial(m)n=[n,\partial(m)],\\
		&[\partial(m),m']_*=\partial(m)*_1 m'-m'*_2\partial(m)=mm'-m'm=[m,m'],
		\end{align*}
		where we use \eqref{XAs1} to prove \eqref{XLie1} and \eqref{XLie2} to prove \eqref{XAs2}.
	\end{proof}
\end{prop}

\begin{remark}
	With the last property we can see that the examples given for the associative algebras case, $(M,M,(*,*),\Id_M)$,
and for the Lie case for the associate Lie algebra $M^{\mathcal{L}}$, $(M^\mathcal{L},M^{\mathcal{L}},[-,-],\Id_{M^\mathcal{L}})$, are related.
\end{remark}

	We will denote by $\textbf{\textit{X}}(\textbf{\textit{LieAlg}}_K)$ the category of crossed modules of Lie $K$-algebras and its homomorphisms.

\begin{remark}
	The previous proposition give us a functor \begin{center}
		$(-)_\mathcal{X}^\mathcal{L}\colon \textbf{\textit{X}}(\textbf{\textit{AssAlg}}_K)\xrightarrow{}\textbf{\textit{X}}(\textbf{\textit{LieAlg}}_K)$.
	\end{center}
\end{remark}

We have the next proposition which relates the categorical Lie $K$-algebras and the categorical associative $K$-algebras.

\begin{prop}
	If $(C_1,C_0,s,t,e,k)$ is a categorical associative $K$-algebra, then $(C_1^\mathcal{L},C_0^\mathcal{L},s,t,e,k)$ is a categorical Lie $K$-algebra.
	\begin{proof}
		Immediate since $(C_1\times_{C_0}C_1)^\mathcal{L}=C_1^{\mathcal{L}}\times_{C_0^{\mathcal{L}}}C_1^{\mathcal{L}}$. They are the same underlying vector space and have the same operation.
	\end{proof}
\end{prop}

\begin{remark}
	The previous proposition give us a functor \begin{center}
		$(-)_\mathcal{C}^\mathcal{L}\colon \textbf{\textit{ICat}}(\textbf{\textit{AssAlg}}_K)\xrightarrow{}\textbf{\textit{ICat}}(\textbf{\textit{LieAlg}}_K)$.
	\end{center}
\end{remark}

\begin{remark}
	As in the case of groups and associative $K$-algebras, the categories $\textbf{\textit{ICat}}(\textbf{\textit{LieAlg}}_K)$ and $\textbf{\textit{X}}(\textbf{\textit{LieAlg}}_K)$ are equivalent (see~\cite{ThRa}).
	
	It is immediate to check that the equivalence functors commute with the functors $(-)_\mathcal{X}^\mathcal{L}$ and 	$(-)_\mathcal{C}^\mathcal{L}$.
 The only problem can be the semidirect product, but it will be solved in the following proposition.
\end{remark}

\begin{defi}	
	Let $M$ and $N$ be two associative $K$-algebras and $\cdot$ an associative action of $N$ on $M$. We define its \emph{semidirect product} as the $K$-vector space $M\times N$ with the following operation:
	\begin{equation*}
	(m,n)(m',n')=(mm'+n\cdot_1 m'+m\cdot_2 n',nn').
	\end{equation*}
	
	Let $M$ and $N$ be two Lie $K$-algebras and $\cdot$ a Lie action of $N$ on $M$. We define its \emph{semidirect product} as the $K$-vector space $M\times N$ with the following bracket:
	\begin{equation*}
	[(m,n),(m',n')]=([m,m']+n\cdot m'-n'\cdot m,[n,n']).
	\end{equation*}
	
	In both cases we will denote the semidirect product by $M\rtimes N$.
\end{defi}

\begin{prop}
Let $M$ and $N$ be associative $K$-algebras. Then, if $*$ is an associative action of $N$ on $M$ (then $[-,-]_*$ is a Lie action of
 $N^\mathcal{L}$ in $M^\mathcal{L}$) we have that $(A\rtimes B)^\mathcal{L}=A^\mathcal{L}\rtimes B^\mathcal{L}$.
\begin{proof}
Since the underlying vector space is the same, we only need to prove that the bracket is the same.
\begin{align*}
&(m,n)(m',n')-(m',n')(m,n)\\
&=(mm'+n*_1m'+m*_2n',nn')-(m'm+n'*_1m+m'*_2n,n'n)\\
&=(mm'+n*_1m'+m*_2n'-m'm-n'*_1 m-m'*_2 n,nn'-n'n)\\
&=([m,m']+[n,m']_*-[n',m]_*,[n,n']),
\end{align*}
 where $(m,n),(m',n')\in M\times N$.
\end{proof}
\end{prop}

\subsection{Braided Semigroupal Category}\hfill
	
	A bifunctor is a functor whose source category is a product category.
	
	Let $F\colon \textbf{\textit{C}}\times \textbf{\textit{D}}\xrightarrow{}\textbf{\textit{E}}$ be a bifunctor.
  For $A\in \Ob(\textbf{\textit{C}})$ and $B\in\Ob(\textbf{\textit{D}})$, we denote by ${}_AF$ and $F_B$ the functors:
	\begin{align*}
	{}_AF\colon \textbf{\textit{D}}\xrightarrow{}\textbf{\textit{E}},\ {}_AF(D\xrightarrow{f}D')&=F(A,D)\xrightarrow{F(\Id_A,f)}F(A,D'),\\
	F_B\colon \textbf{\textit{C}}\xrightarrow{}\textbf{\textit{E}},\ F_B(C\xrightarrow{g}C')&=F(C,B)\xrightarrow{F(g,\Id_B)}F(C',B).
	\end{align*}

Crane and Yetter defined in~\cite{CrYet98} the notion of semigroupal category.

\begin{defi}
	A \emph{semigroupal category} is a triple $\textbf{\textit{C}}=(\textbf{\textit{C}},\otimes,a)$
where \textbf{\textit{C}} is a category, $\otimes\colon \textbf{\textit{C}}\times\textbf{\textit{C}}\xrightarrow{}\textbf{\textit{C}}$ is a bifunctor and
 $a\colon \otimes \circ(\otimes\times\Id_{\textbf{\textit{C}}})\xrightarrow{}\otimes \circ (\Id_{\textbf{\textit{C}}}\times \otimes)$
 is a natural isomorphism called associativity, which verifies the following associative coherence diagram:
	\[
	\begin{tikzcd}
	((X\otimes Y)\otimes Z)\otimes W \arrow[d,"a_{X,Y,Z}\otimes \Id_W"] \arrow[rr,"a_{X\otimes Y,Z,W}"] && (X\otimes Y) \otimes (Z \otimes W) \arrow[dd,"a_{X,Y,Z\otimes W}"]\\
	(X\otimes (Y\otimes Z))\otimes W \arrow[d,"a_{X,Y\otimes Z,W}"]\\
	X\otimes ((Y\otimes Z)\otimes W) \arrow[rr,"\Id_X\otimes a_{Y,Z,W}"] && X\otimes (Y\otimes (Z\otimes W)).
	\end{tikzcd}
	\]
	We will say that a semigroupal category is \emph{strict} if the isomorphism $a$ is the identity morphism. In this case we have that $(X\otimes Y)\otimes Z=X\otimes (Y\otimes Z)$.
\end{defi}

The definition of monoidal category was given in the works~\cite{Ben63,MLane63}.

\begin{defi}
	A \emph{monoidal category} is a $6$-tuple $\textbf{\textit{C}}=(\textbf{\textit{C}},\otimes,a,I,l,r)$ where $(\textbf{\textit{C}},\otimes ,a)$
 is a semigroupal category, $I\in\Ob(\textbf{\textit{C}})$ and $l\colon {}_I\otimes\xrightarrow{}\Id_{\textbf{\textit{C}}}$, $r\colon \otimes_I\xrightarrow{}\Id_{\textbf{\textit{C}}}$
  are natural isomorphisms called, respectively, left unit and right unit constraints, which also verify, for $X,Y\in \Ob(\textbf{\textit{C}})$, the unit coherence diagram:
	\[
	\begin{tikzcd}
	(X\otimes I)\otimes Y\arrow[rd,"r_X\otimes\Id_Y"']\arrow[rr,"a_{X,I,Y}"] && X\otimes (I\otimes Y)\arrow[dl,"\Id_X\otimes l_Y"]\\
	& X\otimes Y.
	\end{tikzcd}
	\]
		We will say that a monoidal category is \emph{strict} if the isomorphisms $a$, $l$ and $r$ are the identity morphisms.
In this case we have that $(X\otimes Y)\otimes Z=X\otimes (Y\otimes Z)$, $X\otimes I=X=I\otimes X$.
\end{defi}	

	If \textbf{\textbf{\textit{C}}} is a category we have the functor
 $T\colon \textbf{\textbf{\textit{C}}}\times \textbf{\textbf{\textit{C}}}\xrightarrow{} \textbf{\textbf{\textit{C}}}\times \textbf{\textbf{\textit{C}}}$ given by the expression
	\begin{equation*}
	T((A,B)\xrightarrow{(f,g)}(A',B'))\coloneqq (B,A)\xrightarrow{(g,f)}(B',A').
	\end{equation*}
	It is easy to check that it is a functorial isomorphism with $T \circ T=\Id_{\textbf{\textit{C}}}$.

In~\cite{JSBM}, Joyal and Street introduced the notion of braided monoidal category.

\begin{defi}
	A \emph{braiding on a monoidal category} \textbf{\textit{C}} is a natural isomorphism $\tau\colon \otimes \xrightarrow{} \otimes \circ T$
such as for any objects $X,Y,Z$ in \textbf{\textit{C}} the following diagrams (associative coherence) commute:

	\[\begin{tikzcd}
	(X\otimes Y)\otimes Z\arrow[r,"\tau_{X\otimes Y,Z}"]\arrow[d,"a_{X,Y,Z}"] & Z\otimes (X\otimes Y)\\
	X\otimes (Y\otimes Z)\arrow[d,"\Id_X\otimes \tau_{Y,Z}"] & (Z\otimes X)\otimes Y\arrow[u,"a_{Z,X,Y}"]\\
	X\otimes (Z\otimes Y)\arrow[r,"a_{X,Z,Y}^{-1}"] & (X\otimes Z)\otimes Y\arrow[u,"\tau_{X,Z}\otimes \Id_Y"],
	\end{tikzcd}\begin{tikzcd}
	X\otimes (Y\otimes Z)\arrow[r,"\tau_{X,Y\otimes Z}"]\arrow[d,"a^{-1}_{X,Y,Z}"] & (Y\otimes Z)\otimes X\\
	(X\otimes Y)\otimes Z\arrow[d,"\tau_{X,Y}\otimes \Id_Z"] & Y\otimes(Z\otimes X)\arrow[u,"a^{-1}_{Y,Z,X}"]\\
	(Y\otimes X)\otimes Z\arrow[r,"a_{Y,X,Z}"] & Y\otimes (X\otimes Z)\arrow[u,"\Id_Y\otimes\tau_{X,Z}"],
	\end{tikzcd}\]
\end{defi}

We can define the concept of braided semigroupal category just by imitating the definition given for the case of monoidal categories.

\begin{defi}
	A\emph{ braiding on a semigroupal category} \textbf{\textit{C}} is a natural isomorphism $\tau\colon \otimes \xrightarrow{} \otimes \circ T$
which verifies the two associative coherence diagrams given in the definition of braiding for monoidal categories.
\end{defi}
\section{Braiding for categorical groups and crossed modules of groups}\label{S:braidgrp}

We have the following property whose proof can be seen in~\cite{ThRa} for a general case.
\begin{lem}\label{Lemacomposition}
	We will suppose that $(C_1,C_0,s,t,e,k)$ is a categorical associative, Lie or Leibniz $K$-algebra or a categorical group (where the operation in $C_1$ is denoted by ``$+$'').
 Then, if $(x,y)\in C_1\times_{C_0}C_1$, the following rule for the composition is true:
	\[k((x,y))=x-e(t(x))+y=x-e(s(y))+y.\]
\end{lem}

\begin{lem}
	In the categories of categorical associative, Lie or Leibniz $K$-algebras and in the category of categorical groups all internal morphisms $f\in C_1$ are internal isomorphisms.
 That is, there exists $f'\in C_1$ such that $k((f,f'))=e(s(f))$ and $k((f',f))=e(t(f))$.
\end{lem}

It is the same to give a strict monoidal category over a small category (internal category in the case of sets) than to give a categorical monoid.
The correspondence is given by taking the product of the monoids $C_0$ and $C_1$ as the $\otimes$ product.  The fact that the morphisms $s$, $t$, $e$ and $k$
are homomorphisms of monoids is equivalent to the functoriality of $\otimes$. The unit $I$ is given by the monoid unit $1_{C_0}$, and the unit in $C_1$ is $e(I)$.

Using this idea one can define what is a braiding for a
categorical monoid, but if we think in the case of groups (a little more restrictive) we have that all internal morphisms are isomorphisms. Using this, it is only needed to take a family of internal morphisms.

The definition of braiding on a categorical group was introduced by Joyal and Street in~\cite{JSBM} and~\cite{JS93}.

\begin{defi}
	Let $\mathcal{C}=(C_1,C_0,s,t,e,k)$ be a categorical group. A \emph{braiding in $\mathcal{C}$} is a map $\tau\colon C_0\times C_0\xrightarrow{}C_1$, $(a,b)\mapsto \tau_{a,b}$,
 such that the following properties are verified:
	\begin{equation}\label{GrT1}
	\tau_{a,b}\colon ab\xrightarrow{ }ba,\tag{GrT1}
	\end{equation}
	\begin{equation}\label{GrT2}
	{\begin{tikzcd}
	{s(x)s(y)}\arrow[d,"{\tau_{s(x),s(y)}}"]\arrow[r,"{xy}"]& {t(x)t(y)}\arrow[d,"{\tau_{t(x),t(y)}}"]\\
	{s(y)s(x)}\arrow[r,"{yx}"]& {t(y)t(x)}.
	\end{tikzcd}}\tag{GrT2}
	\end{equation}
	\begin{align}
	&\tau_{ab,c}=(\tau_{a,c}e(b))\circ (e(a)\tau_{b,c}).\tag{GrT3}\label{GrT3}\\
	&\tau_{a,bc}=(e(b)\tau_{a,c})\circ (\tau_{a,b}e(c)).\tag{GrT4}\label{GrT4}
	\end{align}	
	for $a,b,c\in C_0$, $x,y\in C_1$.
	
	We will say that $(C_1,C_0,s,t,e,k,\tau)$ is a \emph{braided categorical group}.
\end{defi}

\begin{remark}
	One can see that \eqref{GrT2} is the naturality and \eqref{GrT3}, \eqref{GrT4}
are the coherence diagrams.
\end{remark}

\begin{defi}
	A \emph{braided internal functor between two braided categorical groups}, whose braidings are $\tau$ and $\tau'$,
is an internal functor $(F_1,F_0)$ between the internal categories verifying $F_1(\tau_{a,b})=\tau'_{F_0(a),F_0(b)}$ for $a,b\in C_0$.
	
	We denote the category of braided categorical groups and braided internal functors between them as $\textbf{\textit{BICat}}(\textbf{\textit{Grp}})$.
\end{defi}

The definition of braiding in crossed modules of groups was given by Conduch\'e in~\cite[Equalities (2.12)]{Condu} and, although in this case the action is superfluous, it can be recovered as he says previously
as $m\cdot l=l\{\partial(l)^{-1},m\}$.
We will take this action into account, so we double one of the equalities; and we use the last two of equalities (2.11) of \cite{Condu} rather that the last two of (2.12).
 Although there is not problem because Conduch\'e proves in that cite that are equivalent.

\begin{defi}
	Let $(G,H,\cdot, \partial)$ be a crossed module of groups. A \emph{braiding} (or \emph{Peiffer lifting}) of that crossed module is a map $\{-,-\}\colon H\times H\xrightarrow{}G$ which verifies:
	\begin{align}
	\partial\{h,h'\}&=[h,h'],\tag{BGr1}\label{BGr1}\\
	\{\partial g,\partial g'\}&=[g,g'],\tag{BGr2}\label{BGr2}\\
	\{\partial g,h\}&=g(h\cdot g^{-1}),\tag{BGr3}\label{BGr3}\\
	\{h,\partial g\}&=(h\cdot g)g^{-1},\tag{BGr4}\label{BGr4}\\
	\{h,h'h''\}&=\{h,h'\}(h'\cdot\{h,h''\}),\tag{BGr5}\label{BGr5}\\
	\{hh',h''\}&=(h\cdot\{h',h''\})\{h,h''\},\tag{BGr6}\label{BGr6}
	\end{align}
	for $g,g'\in G$, $h,h',h''\in H$, where $[g,g']=gg'g^{-1}g'^{-1}$.
	
	We say that $(G,H,\cdot, \partial,\{-,-\})$ is a \emph{braided crossed module of groups}.
\end{defi}

\begin{example}
	It is easy to check that the
	commutator $[-,-]$ is a braiding on $(G,G,\Conj,\Id_G)$.
\end{example}

\begin{defi}
	 $(G,H,\cdot,\partial,\{-,-\})\xrightarrow{(f_1,f_2)}(G',H',*,\partial',\{-,-\}')$ is an \emph{homomorphism of braided crossed modules of groups}
 if it is an homomorphism of crossed modules of groups such that $f_1(\{h,h'\})=\{f_2(h),f_2(h')\}'$ for $h,h'\in H$.
	
	We denote the category of  braided crossed modules of groups and its homomorphisms as $\textbf{\textit{BX}}(\textbf{\textit{Grp}})$.
\end{defi}

\begin{remark}
	We can see in~\cite{JSBM,JS93} and~\cite{Garzon&Miranda} that the categories $\textbf{\textit{BICat}}(\textbf{\textit{Grp}})$ and $\textbf{\textit{BX}}(\textbf{\textit{Grp}})$ are equivalent.
\end{remark}

\section{Braiding for categorical associative algebras and crossed modules of associative algebras}\label{S:braidassalg}

In this section we will introduce a definition of braiding for categorical associative $K$-algebras.

	As categorical monoids can be seen as strict monoidal internal categories in \textbf{\textit{Set}},
 we can think that a strict semigroupal category over an internal category in \textbf{\textit{Vect}}$_K$ is really a categorical associative $K$-algebra.
  By the same reasoning we have that we can identify the $\otimes$ product with the second operation and the functoriality is recovered in the same way as the case of groups.
   The $K$-bilinearity of the product is given by the fact that $\otimes$ is a internal bifunctor in \textbf{\textit{Vect}}$_K$ and, this means,
   it is a functor between the respective small categories and each component (fixing an object on left or right) is an internal functor, that means, is linear in internal objects and internal morphisms.
	
	With this in mind we can introduce a braiding on categorical associative $K$-algebras, imitating the braiding for semigroupal categories.
	
	We will add that $\tau\colon C_0\times C_0\xrightarrow{} C_1$ is $K$-bilinear, but it is obvious if we have the idea that for an internal object $A\in C_0$
we must have the morphisms in \textbf{\textit{Vect}}$_K$ $\tau_{A,-},\tau_{-,A}\colon C_0\xrightarrow{}C_1$ defined by $\tau_{A,-}(B)=\tau_{A,B}$ and $\tau_{-,A}(B)=\tau_{B,A}$.
	
	We must recall that we show that in the internal categories in which we will work all internal morphisms are internal isomorphisms.

\begin{defi}
	Let $\mathcal{C}=(C_1,C_0,s,t,e,k)$ be a categorical associative $K$-algebra.
	
	A\emph{ braiding on $\mathcal{C}$} is a $K$-bilinear map $\tau\colon C_0\times C_0\xrightarrow{} C_1$, $(a,b)\mapsto \tau_{a,b}$, verifying the following properties:
	\begin{equation}\label{AsT1}
	\tau_{a,b}\colon ab\xrightarrow{ }ba,\tag{AsT1}
	\end{equation}
	\begin{equation}\label{AsT2}
	\begin{tikzcd}
	{s(x)s(y)}\arrow[d,"{\tau_{s(x),s(y)}}"]\arrow[r,"{xy}"]& {t(x)t(y)}\arrow[d,"{\tau_{t(x),t(y)}}"]\\
	{s(y)s(x)}\arrow[r,"{yx}"]& {t(y)t(x)},
	\end{tikzcd}\tag{AsT2}
	\end{equation}
	\begin{align}
	\tau_{ab,c}&=(\tau_{a,c}e(b))\circ (e(a)\tau_{b,c}),\tag{AsT3}\label{AsT3}\\
	\tau_{a,bc}&=(e(b)\tau_{a,c})\circ (\tau_{a,b}e(c)),\tag{AsT4}\label{AsT4}
	\end{align}
	for $a,b,c\in C_0$, $x,y\in C_1$.
	
	We will say that $(C_1,C_0,s,t,e,k,\tau)$ is a \emph{braided categorical associative $K$-algebra}.
	\end{defi}

\begin{remark}
	As in the case of groups, \eqref{AsT2} is the naturality and \eqref{AsT3}, \eqref{AsT4} are the coherence diagrams.
\end{remark}

\begin{defi}
	A \emph{braided internal functor between two braided categorical associative $K$-algebras},
whose braidings are $\tau$ and $\tau'$, is an internal functor $(F_1,F_0)$ such that $F_1(\tau_{a,b})=\tau'_{F_0(a),F_0(b)}$ for $a,b\in C_0$.
\end{defi}

We denote the category of braided categorical associative $K$-algebras and braided internal functors between them as $\textbf{\textit{BICat}}(\textbf{\textit{AssAlg}}_K)$.

We will introduce the notion of braiding for crossed modules of associative algebras looking for an equivalence between braided crossed modules and braided internal categories of associative algebras. as it happens in the groups case,
 The definition we got is the following one.

\begin{defi}
	Let $(M\xrightarrow{\partial}N,*)$ be a crossed module of associative $K$-algebras, where $*=(*_1,*_2)$.
 Then a \emph{braiding} on the crossed module is a $K$-bilinear map $\{-,-\}\colon N\times N\xrightarrow{}M$ verifying:
	\begin{align}
	\partial\{n,n'\}&=[n,n'],\tag{BAs1}\label{BAs1}\\
	\{\partial m, \partial m' \}&=[m,m'],\tag{BAs2}\label{BAs2}\\
	\{\partial m, n \}&=-[n,m]_*, \tag{BAs3}\label{BAs3}\\
	\{n,\partial m \}&=[n,m]_*, \tag{BAs4}\label{BAs4}\\
	\{n,n'n''\}&=n'*_1\{n,n''\}+\{n,n'\}*_2n'',\tag{BAs5}\label{BAs5}\\
	\{nn',n''\}&=n*_1\{n',n''\}+\{n,n''\}*_2n',\tag{BAs6}\label{BAs6}
	\end{align}
	for $m,m'\in M$, $n,n',n''\in N$.
	
	We denote $[n,m]_*=n*_1m-m*_2n$ and $[x,y]=xy-yx$.
	
	If $\{-,-\}$ is a braiding we will say that $(M\xrightarrow{\partial}N,*,\{-,-\})$ is a \emph{braided crossed module of associative $K$-algebras}.
\end{defi}

\begin{example}
	We have that the commutator $[-,-]$ is a braiding on the crossed module $(M,M,(*,*),\Id_M)$.
\end{example}

\begin{defi}
	An \emph{homomorphism of braided crossed modules of associative $K$-algebras} $(M,N,\cdot,\partial,\{-,-\})\xrightarrow{(f_1,f_2)}(M',N',*,\partial',\{-,-\}')$
 is an homomorphism of crossed modules of associative $K$-algebras such that $f_1(\{n,n'\})=\{f_2(n),f_2(n')\}'$ for $n,n'\in N$.
\end{defi}

We denote the category of  braided crossed modules of associative $K$-algebras and its homomorphisms as $\textbf{\textit{BX}}(\textbf{\textit{AssAlg}}_K)$.

\begin{prop}
	Let $\mathcal{X}=(M,N,(*_1,*_2),\partial,\{-,-\})$ be a braided crossed module of associative $K$-algebras.
	
	Then $\mathcal{C}_\mathcal{X} \coloneqq (M\rtimes N,N,\bar{s},\bar{t},\bar{e},\bar{k},\bar{\tau})$ is a braided categorical associative $K$-algebra where:
	\begin{itemize}
		\item $\bar{s}\colon M\rtimes N\xrightarrow{} N$, $\bar{s}((m,n))=n$,
		\item $\bar{t}\colon M\rtimes N\xrightarrow{} N$, $\bar{t}((m,n))=\partial m+n$,
		\item $\bar{e}\colon N\xrightarrow{} M\rtimes N$, $\bar{e}(n)=(0,n)$,
		\item $\bar{k}\colon (M\times N)\times_N (M\rtimes N)\xrightarrow{}M\rtimes N$, where the source is the pullback of $\bar{t}$ along $\bar{s}$, defined as $k(((m,n),(m',\partial m+n)))=(m+m',n)$.
		\item $\bar{\tau}\colon N\times N\xrightarrow{}M\rtimes N$, $\bar{\tau}_{n,n'}=(-\{n,n'\},nn')$.
	\end{itemize}
	\begin{proof}
		Other than the braiding, it has already been proven, as can be seen  in~\cite{ThRa}, that $(M\rtimes N,N,\bar{s},\bar{t},\bar{e},\bar{k})$ is a categorical associative $K$-algebra.
 We only need to check the braiding axioms for this internal category.
		
		We will start with \eqref{AsT1}. Let us take $n,n'\in N$.
		\begin{align*}
		\bar{s}(\bar{\tau}_n,n')=\bar{s}((-\{n,n'\},nn'))&=nn',\\
		\bar{t}(\bar{\tau}_{n,n'})=\bar{t}((-\{n,n'\},nn'))&=-\partial\{n,n'\}+nn'\\
		&=-[n,n']+nn'=n'n,
		\end{align*}
		where we use \eqref{BAs1}.
		
		We will prove now \eqref{AsT2}. We will take $x=(m,n),y=(m',n')\in M\rtimes N$.
		
		We need to show that $\tau_{t(x),t(y)}\circ xy=yx\circ \tau_{s(x),s(y)}$.
		\begin{align*}
		&\tau_{t(x),t(y)}\circ xy\\
		&=\bar{k}(((m,n)(m',n'),(-\{\bar{t}((m,n)),\bar{t}((m',n'))\},\bar{t}((m,n))\bar{t}((m',n')))))\\
		&=\bar{k}(((m,n)(m',n'),(-\{\partial m+n,\partial m'+n'\},(\partial m+n)(\partial m'+n'))))\\
		&=\bar{k}(((mm'+n*_1 m'+m*_2 n',nn'),(-\{\partial m+n,\partial m'+n'\},\\
		&\qquad\qquad\qquad\qquad\qquad\qquad\qquad\qquad\qquad\qquad(\partial m+n)(\partial m'+n'))))\\
		&=(mm'+n*_1 m'+m*_2n'-\{\partial m+n,\partial m'+n'\},nn')  \\
		&=(mm'+n*_1 m'+m*_2n'-\{\partial m,\partial m'\}-\{\partial m,n'\}-\{n,\partial m'\} -\{n,n'\},nn')\\			
		&=(mm'+n*_1 m'+m*_2n'-[m,m']+[n',m]_*-[n,m']_* -\{n,n'\},nn')\\
		&=(m'm+m'*_2 n+n'*_1 m-\{n,n'\},nn'),
		\end{align*}
		where we use \eqref{BAs2}, \eqref{BAs3} and \eqref{BAs4} in the sixth equality. In the other way,
		\begin{align*}
		& yx\circ \tau_{s(x),s(y)}\\
		&=\bar{k}(((-\{\bar{s}((m,n)),\bar{s}((m,n'))\},\bar{s}((m,n))\bar{s}((m',n'))),(m',n')(m,n)))\\
		&=\bar{k}(((-\{n,n'\},nn'),(m',n')(m,n)))\\	
		&=\bar{k}(((-\{n,n'\},nn'),(m'm+n'*_1 m+m'*_2 n,n'n)))\\
		&=(-\{n,n'\}+m'm+n'*_1 m+m'*_2 n,nn'),
		\end{align*}
		where we can see the equality.
		
		We will verify \eqref{AsT3} below. If $n,n',n''\in N$, then
		\begin{align*}
		&(\bar{\tau}_{n,n''}\bar{e}(n'))\circ (\bar{e}(n)\bar{\tau}_{n',n''})=\bar{k}((\bar{e}(n)\bar{\tau}_{n',n''},\bar{\tau}_{n,n''}\bar{e}(n')))\\
		&=\bar{k}((0,n)(-\{n',n''\},n'n''),(-\{n,n''\},nn'')(0,n'))\\
		&=\bar{k}(-n*_1\{n',n''\},n(n'n''),(-\{n,n''\}*_2 n',(nn'')n'))\\
		&=(-n*_1\{n',n''\}-\{n,n''\}*_2n',n(n'n''))=(-\{nn',n''\},(nn')n'')=\bar{\tau}_{nn',n''},
		\end{align*}
		where we use \eqref{BAs6} and associativity.
		
		Finally we will show that \eqref{AsT4} is verified.
		If $n,n',n''\in N$, then
		\begin{align*}
		&(\bar{e}(n')\bar{\tau}_{n,n''})\circ (\bar{\tau}_{n,n'}\bar{e}(n''))=\bar{k}(\tau_{n,n'}e(n''),e(n')\tau_{n,n''})\\
		&=\bar{k}((-\{n,n'\},nn')(0,n''),(0,n')(-\{n,n''\},nn''))\\
		&=\bar{k}((-\{n,n'\}*_2 n'',(nn')n''),(-n'*_1\{n,n''\},n'(nn'')))\\
		&=(-n'*_1\{n,n''\}-\{n,n'\}*_2 n'',(nn')n'')=(-\{n,n'n''\},n(n'n''))=\bar{\tau}_{n,n'n''},
		\end{align*}
		where we use \eqref{BAs5} along with the associativity in the second equality.
		
		Thus the braiding axioms are verified for the categorical associative $K$-algebra.
	\end{proof}
\end{prop}

\begin{prop}
	We have a functor $\mathcal{C_\mathfrak{A}}\colon \textbf{\textit{BX}} (\textbf{\textit{AssAlg}}_K) \xrightarrow{} \textbf{\textit{BICat}} ( \textbf{\textit{AssAlg}}_K)$ defined as:
	\[\mathcal{C_\mathfrak{A}}(\mathcal{X}\xrightarrow{(f_1,f_2)}\mathcal{X}')=\mathcal{C}_\mathcal{X}\xrightarrow{(f_1\times f_2,f_2)}\mathcal{C}_{\mathcal{X}'}\]
	where $\mathcal{C}_\mathcal{X}$ is defined in the previous proposition.
	\begin{proof}
		We know that the pair $(f_1\times f_2,f_2)$ is an internal functor between the respective internal categories,
 since what we are trying to do is to extend an existing functor (see~\cite{ThRa}) to the braided case. In the same way, as it is an extension, we already know that,
  if it is well defined, it verifies the properties of functor, since the composition and identity are the same as in the categories without braiding.
		
		Because of that, to conclude this proof, is enough to see that $(f_1\times f_2,f_2)$ is a braided internal functor of braided categorical associative $K$-algebras.	
		\begin{align*}
		&(f_1\times f_2)(\bar{\tau}_{n,n'})=(f_1\times f_2)((-\{n,n'\},nn'))=(-f_1(\{n,n'\}),f_2(nn'))\\
		&=(-\{f_2(n),f_2(n')\}',f_2(n)f_2(n'))=\bar{\tau}'_{f_2(n),f_2(n')},
		\end{align*}
		where we use that $(f_1,f_2)$ is an homomorphism of braided crossed modules of associative algebras.
	\end{proof}
\end{prop}

\begin{prop}
	Let $\mathcal{C}=(C_1,C_0,s,t,e,k,\tau)$ be a braided categorical associative $K$-algebra.
	
	Then $\mathcal{X}_\mathcal{C} \coloneqq (\ker(s),C_0,({}^{e}*,*^e),\partial_t,\{-,-\}_\tau)$ is a braided crossed module of associative $K$-algebras, where:
	\begin{itemize}
		\item ${}^{e}*\colon C_0\times \ker(s)\xrightarrow{}\ker(s)$, $a\ {{}^{e}*} \ x \coloneqq e(a)x$,
		\item $*^e\colon \ker(s)\times C_0\xrightarrow{}\ker(s)$, $x*^e a \coloneqq xe(a)$,
		\item $\partial_t \coloneqq t|_{\ker(s)}$,
		\item $\{-,-\}_\tau\colon C_0\times C_0\xrightarrow{}\ker(s)$, $\{a,b\}_\tau \coloneqq e(ab)-\tau_{a,b}$.
	\end{itemize}
	
	\begin{proof}
		It is proven (see~\cite{ThRa}) that, under these hypotheses, $(\ker(s),C_0,(*^e,{}^{e}*),\partial_t)$ is a crossed module of associative $K$-algebras.
		
		For that, is enough to show that $\{-,-\}_\tau$ is a braiding on that crossed module.
		
		First we see that it is well defined,
		that means $\{a,b\}_\tau\in\ker(s)$ for $a,b\in C_0$.
		\begin{align*}
		s(\{a,b\}_\tau)=s(e(ab)-\tau_{a,b})=ab-ab=0,
		\end{align*}
		where we use \eqref{AsT1}.
		
		Since they are well defined, we can see if they verify the properties.
		
		To start we will check \eqref{BAs1}. If $a,b\in C_0$, then
		\begin{align*}
		\partial_{t}\{a,b\}_\tau=t(e(ab)-\tau_{a,b})=ab-ba=[a,b],
		\end{align*}
		where we use \eqref{AsT1}.
		
		We will see if it is verified \eqref{BAs2}. If $x,y\in \ker(s)$, then
		\begin{align*}
		\{\partial_tx,\partial_ty\}_\tau=e(\partial_tx\partial_ty)-\tau_{\partial_t x,\partial_t y}
		=e(t(x)t(y))-\tau_{t(x),t(y)}.
		\end{align*}
		We need to show that $e(t(x)t(y))-\tau_{t (x),t (y)}=[x,y]$.
		
		By the axiom \eqref{AsT2} we know the equality
		\[k((xy,\tau_{t(x),t(y)}))=k((\tau_{s(x),s(y)},yx)).\]
		As $x\in \ker(s)$, we have that $s(x)=0$ (in the same way $y$), and $\tau_{s(x),s(y)}=0$ by $K$-bilinearity. We have then that
		\[k((\tau_{s(x),s(y)},yx))=k((0,yx)),\]
		and therefore the equality
		\[k((xy,\tau_{t(x),t(y)}))=k((0,yx)).\]
		Using now the $K$-linearity of $k$ in the previous expression, we obtain
		\[0=k((xy,\tau_{t(x),t(y)}-yx)).\]
		
		Since $t(\tau_{t(x),t(y)}-yx)=t(y)t(x)-t(y)t(x)=0=s(e(0))$ we can talk about $k((\tau_{t(x),t(y)}-yx,e(0)))$.
		Further $k((\tau_{t(x),t(y)}-yx,e(0)))=\tau_{t(x),t(y)}-yx$ by the internal category axioms.
		
		Adding both equalities and using the $K$-linearity of $k$, we get
		\[k((xy+\tau_{t(x),t(y)}-yx,\tau_{t(x),t(y)}-yx))=\tau_{t(x),t(y)}-yx.\]
		Therefore, by grouping, we have
		\[k(([x,y]+\tau_{t(x),t(y)},\tau_{t(x),t(y)}-yx))=\tau_{t(x),t(y)}-yx.\]
		As $s(\tau_{t(x),t(y)}-yx)=t(x)t(y)+0=t(x)t(y)$ (we use that $x$ or $y$ are in $\ker(s)$)
it makes sense to talk about the composition $k((e(t(x)t(y)),\tau_{t(x),t(y)}-yx))$, which is equal to $\tau_{t(x),t(y)}-yx$.
		
		Subtracting both equalities and using the $K$-linearity of $k$, we obtain
		\[k(([x,y]+\tau_{t(x),t(y)}-e(t(x)t(y)),0))=0.\]
		
		Again, using the properties for internal categories, we have
		\begin{align*}
		0&=k(([x,y]+\tau_{t(x),t(y)}-e(t(x)t(y)),0))\\&=k(([x,y]+\tau_{t(x),t(y)}-e(t(x)t(y)),e(0)))
		\\&=[x,y]+\tau_{t(x),t(y)}-e(t(x)t(y)),
		\end{align*}
		which gives us the required equality.
		
		As an observation to the above, in the part of the proof that is related to $x,y\in\ker(s) $, it is sufficient that one of the two is in that kernel.
		Therefore, by repeating the proof using this, we have the following equalities for $x \in \ker(s) $ and $y \in C_1 $:
		\begin{align*}
		e(t(x)t(y))-\tau_{t (x),t (y)}=[x,y], \quad e(t(y)t(x))-\tau_{t (y),t (x)}=[y,x].
		\end{align*}
		
		With these equalities, we will prove \eqref{BAs3} and \eqref{BAs4}.
		
		Let $a\in C_0$ and $x\in \ker(s)$. Then
		\begin{align*}
		\{\partial_t x, a\}_\tau&=e(t(x)t(e(a)))-\tau_{t (x),t(e(a))}=[x,e(a)]=xe(a)-e(a)x=x*^e a-a\ {{}^{e}*}\ x,\\
		\{ a,\partial_t x\}_\tau&=e(t(e(a))t(x))-\tau_{t (e(a)),t(x)}=[e(a),x]=e(a)x-xe(a)=a\ {{}^{e}*} \ x-x*^e a.
		\end{align*}
		
		We will see now the last conditions, starting with \eqref{BAs5}. Let $a,b,c\in C_0$.
		\begin{align*}
		&\{a,bc\}_\tau=e(a(bc))-\tau_{a,bc}=e(a(bc))-((e(b)\tau_{a,c})\circ (\tau_{a,b}e(c)))\\
		&=e(a(bc))-e(b)\tau_{a,c}-\tau_{a,b}e(c)+e(t(\tau_{a,b}e(c)))\\
		&=e((ab)c)-e(b)\tau_{a,c}-\tau_{a,b}e(c)+e((ba)c)\\
		&=e(b)e(ac)-e(b)\tau_{a,c}+e(ab)e(c)-\tau_{a,b}e(c)\\
		&=e(b)\{a,c\}_\tau+\{a,b\}_\tau e(c)=b\ {{}^{e}*}\ \{a,c\}_\tau+\{a,b\}_\tau*^e c,
		\end{align*}
		where we use \eqref{AsT4}, Lemma~\ref{Lemacomposition} and the associativity.
		
		To conclude we will check \eqref{BAs6}.
		\begin{align*}
		&\{ab,c\}_\tau=e((ab)c)-\tau_{ab,c}=e((ab)c)-((\tau_{a,c}e(b))\circ (e(a)\tau_{b,c}))\\
		&=e((ab)c)-\tau_{a,c}e(b)-e(a)\tau_{b,c}+e(t(e(a)\tau_{b,c}))\\
		&=e(a(bc))-\tau_{a,c}e(b)-e(a)\tau_{b,c}+e(a(cb))\\
		&=e(a)e(bc)-e(a)\tau_{b,c}+e(ac)e(b)-\tau_{a,c}e(b)\\
		&=e(a)\{b,c\}_\tau + \{a,c\}_\tau e(b)=a\ {{}^{e}*}\ \{b,c\}_\tau+\{a,c\}_\tau *^e b,
		\end{align*}
		where we use \eqref{AsT3}, Lemma~\ref{Lemacomposition} and associativity.
	\end{proof}
\end{prop}

\begin{prop}
	We have a functor $\mathcal{X_\mathfrak{A}}\colon \textbf{\textit{ICat}} (\textbf{\textit{AssAlg}}_K) \xrightarrow{} \textbf{\textit{BX}}(\textbf{\textit{AssAlg}}_K)$ defined as
	\[\mathcal{X_\mathfrak{A}}(\mathcal{C}\xrightarrow{(F_1,F_0)}\mathcal{C}')=\mathcal{X}_\mathcal{C}\xrightarrow{(F_1^s,F_0)}\mathcal{X}_{\mathcal{C}'},\]
	where $\mathcal{X}_\mathcal{C}$ is defined in the previous proposition and $F_1^s\colon \ker(s)\xrightarrow{}\ker(s')$ is defined as $F_1^s(x)=F_1(x)$ for $x\in  \ker(s)$.
	\begin{proof}
		The fact that it is a functor between the categories without braiding  is already shown~\cite{ThRa},
 so we have to see that it can be extended to the braided case. For this, we have to verify the axioms of the homomorphisms of braided crossed modules of associative $K$-algebras.
	
		\begin{align*}
		& F^s_1(\{a,b\}_\tau)=F_1(e(ab)-\tau_{a,b})=F_1(e(a,b))-F_1(\tau_{a,b})\\
		&=e'(F_0(ab))-\tau'_{F_0(a),F_0(b)}=e'(F_0(a)F_0(b))-\tau'_{F_0(a),F_0(b)}\\
		&=\{F_0(a),F_0(b)\}_{\tau'}.\qedhere
		\end{align*}
	\end{proof}
\end{prop}

\begin{remark}
	Note that, if $(M,N,(*_1,*_2),\partial,\{-,-\})$ is a braided crossed module of associative $K$-algebras, then $\ker(\bar{s})=\{(m,0)\in M\rtimes N\mid m\in M\}=:(M,0)$,
 where $\bar{s}$ is defined for the functor $\mathcal{C}_\mathfrak{A}$.
\end{remark}

\begin{prop}
	The categories $\textbf{\textit{BX}}(\textbf{\textit{AssAlg}}_K)$ and $\textbf{\textit{ICat}}(\textbf{\textit{AssAlg}}_K)$ are equivalent categories.
	
	Further, the functors $\mathcal{C_\mathfrak{A}}$ and $\mathcal{X_\mathfrak{A}}$ are inverse equivalences,
 where the natural isomorphisms $\Id_{\textbf{\textit{BX}}(\textbf{\textit{AssAlg}}_K)}\stackrel{\alpha^\mathfrak{A}}{\cong}\mathcal{X}_\mathfrak{A}\circ \mathcal{C}_\mathfrak{A}$ and $\Id_{\textbf{\textit{ICat}}(\textbf{\textit{AssAlg}}_K)}\stackrel{\beta^\mathfrak{A}}{\cong}\mathcal{C}_\mathfrak{A}\circ \mathcal{X}_\mathfrak{A}$ are given by:
	
$\bullet$	If $\mathcal{Z}=(M,N,(\cdot_1,\cdot_2),\partial,\{-,-\})$ is a braided crossed module of associative $K$-algebras,
 then $\alpha^\mathfrak{A}_{\mathcal{Z}}=(\alpha^\mathfrak{A}_M,\Id_N)$ with $\alpha^\mathfrak{A}_M\colon M\xrightarrow{} (M,0)$ defined as $\alpha_M(m)=(m,0)$.
	
$\bullet$	If $\mathcal{D}=(C_1,C_0,s,t,e,k,\tau)$ is a braided categorical associative $K$-algebra,
then $\beta^\mathfrak{A}_{\mathcal{D}}=(\beta^\mathfrak{A}_{s},\Id_{C_0})$ with $\beta^\mathfrak{A}_{C_1}\colon C_1\xrightarrow{}\ker(s)\rtimes C_0$ defined as $\beta^\mathfrak{A}_{C_1}(x)=(x-e(s(x)),s(x))$.
	\begin{proof}
		It can be seen in~\cite{ThRa} that they are well defined maps and that they are isomorphisms in the categories without braiding, as well as they are natural isomorphisms.
		
		For that, it is enough to show that they are isomorphisms between braided objects.
		
		Immediately from the definition, as in the crossed module and internal case, they are isomorphisms if they are bijective morphisms,
 since the inverse map verifies the braided axioms for morphisms in their respective categories.
		
		We know that they are bijective maps, since they are isomorphisms between the categories without braiding. For that, we only have to verify that they are, in fact, morphisms.
		
		Let $\mathcal{Z}=(M,N,(\cdot_1,\cdot_2),\partial,\{-,-\})$ a braided crossed module of associative $K$-algebras.
 Let us see that $\alpha^\mathfrak{A}_{\mathcal{Z}}=(\alpha^\mathfrak{A}_M,\Id_N)$ is an homomorphism.
		\begin{align*}
		&\Id_N(\{n,n'\}_{\bar{\tau}})=\{n,n'\}_{\bar{\tau}}=\bar{e}(nn')-\bar{\tau}_{n,n'}=(0,nn')-(-\{n,n'\},nn')\\
		&=(\{n,n'\},0)=\alpha^\mathfrak{A}_{M}(\{n,n'\}), \quad \text{where} \ n,n'\in N.
		\end{align*}
				
		Let $\mathcal{D}=(C_1,C_0,s,t,e,k,\tau)$ be a braided categorical associative $K$-algebra. We will check that $\beta^\mathfrak{A}_{\mathcal{D}}=(\beta^\mathfrak{A}_{s},\Id_{C_0})$ is a morphism.
		
		If $a,b\in C_0$, we have
		\begin{align*}
		&\Id_{C_0}(\bar{\tau}_{a,b})=\bar{\tau}_{a,b}=(-\{a,b\}_\tau,ab)=(\tau_{a,b}-e(ab),ab)\\
		&=(\tau_{a,b}-e(s(\tau_{a,b})),s(\tau_{a,b}))=\beta^\mathfrak{A}_{C_1}(\tau_{a,b}).
		\end{align*}
		
		Therefore, since they are morphisms, we know that these are natural isomorphisms, as we explained previously, and the equivalence of categories is obtained.
	\end{proof}
\end{prop}

\section{Braiding for categorical Lie algebras and crossed modules of Lie algebras}\label{S:braidLiealg}

In this section we will show that the definition, given by Ulualan in~\cite{Ulua} for braided categorical Lie $K$-algebras,
appears naturally from the previous one using the fact that we can transform an associative $K$-algebra $M$ in a Lie $K$-algebra $M^\mathcal{L}$ with bracket $[x,y]=xy-yx$.

Then we will suppose that $K$ is a field of $\car(K)\neq 2$ to change a little the definition of braiding.
Doing this we will obtain the definition given in~\cite{TFM}, where the equivalence is proven for $\car(K)\neq 2$ with the category of braided crossed modules of Lie $K$-algebras which we will talk about later.

The notion of braiding for categorical Lie $K$-algebras was introduced by Ulualan in~\cite{Ulua}, and it is the following one.

\begin{defi}
Let $\mathcal{C}=(C_1,C_0,s,t,e,k)$ be a categorical Lie $K$-algebra.

A\emph{ braiding on $\mathcal{C}$ } its a $K$-bilinear map $\tau\colon C_0\times C_0\xrightarrow{} C_1$, $(a,b)\mapsto \tau_{a,b}$, which verify the following properties:
\begin{equation}\label{LieT1}
\tau_{a,b}\colon [a,b]\xrightarrow{ }[b,a],\tag{LieT1}
\end{equation}
\begin{equation}\label{LieT2}
\begin{tikzcd}		{[s(x),s(y)]}\arrow[d,"{\tau_{s(x),s(y)}}"]\arrow[r,"{[x,y]}"]& {[t(x),t(y)]}\arrow[d,"{\tau_{t(x),t(y)}}"]\\
{[s(y),s(x)]}\arrow[r,"{[y,x]}"]& {[t(y),t(x)]},\tag{LieT2}
\end{tikzcd}
\end{equation}
\begin{align}
\tau_{[a,b],c}&=[\tau_{a,c},e(b)]+[e(a),\tau_{b,c}],\tag{LieB3}\label{LieB3}\\
\tau_{a,[b,c]}&=[e(b),\tau_{a,c}]+ [\tau_{a,b},e(c)],\tag{LieB4}\label{LieB4}
\end{align}
for $a,b,c\in C_0$, $x,y\in C_1$.

We will say that $(C_0,C_1,s,t,e,k,\tau)$ is a \emph{braided categorical Lie $K$-algebra}.
\end{defi}

\begin{remark}
	It can be seen that, in the previous definition, the lack of associativity in Lie $K$-algebras is fixed in \eqref{LieB3} and \eqref{LieB4} using the Jacobi identity in source and target.
\end{remark}

We want to show that the definition for Lie $K$-algebras is well related with the definition for associative $K$-algebras.

\begin{prop}
	If $(C_1,C_0,s,t,e,k,\tau)$ is a braided categorical associative $K$-algebra, then $(C_1^\mathcal{L},C_0^\mathcal{L},s,t,e,k,\tau^{Lie})$ is a braided categorical Lie $K$-algebra, where
	\[\tau^{Lie} \colon C_0^\mathcal{L}\times C_0^\mathcal{L}\xrightarrow{}C_1^\mathcal{L},\qquad \tau^{Lie}_{a,b} \coloneqq \tau_{a,b}-\tau_{b,a}.\]
	\begin{proof}
		It is easy to see that \eqref{AsT1} implies \eqref{LieT1} and  \eqref{AsT2} implies \eqref{LieT2}.

		By using Lemma~\ref{Lemacomposition} we obtain \eqref{LieB3} and \eqref{LieB4} from \eqref{AsT3} and \eqref{AsT4}.
	\end{proof}
\end{prop}

In~\cite{TFM} is given another definition for braided internal category of Lie $K$-algebras to make the equivalence with the braided crossed modules of Lie $K$-algebras.
 That equivalence is proven for a field with $\car(K)\neq 2$, so we will show that, with that fact, the two definitions are equivalent.

The definition given in~\cite{TFM} is the following one.

\begin{defi}
	Let $\mathcal{C}=(C_1,C_0,s,t,e,k)$ be a categorical Lie $K$-algebra.
	
	A \emph{braiding} on $\mathcal{C}$ is a $K$-bilinear map $\tau\colon C_0\times C_0\xrightarrow{} C_1$, $(a,b)\mapsto \tau_{a,b}$, that verifies \eqref{LieT1}, \eqref{LieT2} and the following properties:
	\begin{align}
	\tau_{[a,b],c}&=\tau_{a,[b,c]}-\tau_{b,[a,c]},\tag{LieT3}\label{LieT3}\\
	\tau_{a,[b,c]}&=\tau_{[a,b],c}-\tau_{[a,c],b},\tag{LieT4}\label{LieT4}
	\end{align}
	for $a,b,c\in C_0$, $x,y\in C_1$.
\end{defi}

The two definitions are equivalent in $\car(K)\neq 2$, as can be seen in the following proposition.

\begin{prop}\label{trenzaanticoncor}
	Let $K$ be a field of $\car(K)\neq 2$ and $(C_1,C_0,s,t,e,k)$ a categorical Lie $K$-algebra.
	
	If $\tau\colon C_0\times C_0\xrightarrow{}C_1$ is a $K$-bilinear map verifying \eqref{LieT1} and \eqref{LieT2}, then
	\begin{align*}
	&\tau_{a,[b,c]}=[e(a),\tau_{b,c}] \quad \text{and}\\
	&\tau_{[b,c],a}=[\tau_{b,c},e(a)].
	\end{align*}
	In particular, by the anticommutativity, we have that $\tau_{a,[b,c]}=-\tau_{[b,c],a}$.
	\begin{proof}
		Using \eqref{LieT1} and \eqref{LieT2} the following diagram is commutative:
		\begin{center}
			\begin{tikzcd}
				{[a,[b,c]]}\arrow[d,"{\tau_{a,[b,c]}}"]\arrow[rr,"{[e(a),\tau_{b,c}]}"]& &{[a,[c,b]]}\arrow[d,"{\tau_{a,[c,b]}}"]\\
				{[[b,c],a]}\arrow[rr,"{[\tau_{b,c},e(a)]}"]& &{[[c,b],a]}.
			\end{tikzcd}
		\end{center}
		That is, we have the equality
		\begin{align*}
		k(([e(a),\tau_{b,c}],\tau_{a,[c,b]}))=k((\tau_{a,[b,c]},[\tau_{b,c},e(a)])),
		\end{align*}
		and then
		\begin{align*}
		0&=k(([e(a),\tau_{b,c}],\tau_{a,[c,b]}))-k((\tau_{a,[b,c]},[\tau_{b,c},e(a)]))\\
		&=k(([e(a),\tau_{b,c}]-\tau_{a,[b,c]},\tau_{a,[c,b]}-[\tau_{b,c},e(a)]))\\
		&=k(([e(a),\tau_{b,c}]-\tau_{a,[b,c]},-\tau_{a,[b,c]}+[e(a),\tau_{b,c}])).
		\end{align*}
		
		Using now Lemma~\ref{Lemacomposition}, we have
		\begin{align*}
		0&=[e(a),\tau_{b,c}]-\tau_{a,[b,c]}+(-\tau_{a,[b,c]}+[e(a),\tau_{b,c}])-e(s(-\tau_{a,[b,c]}+[e(a),\tau_{b,c}]))\\
		&=2([e(a),\tau_{b,c}]-\tau_{a,[b,c]})-e(-[a,[b,c]]+[a,[b,c]])=2([e(a),\tau_{b,c}]-\tau_{a,[b,c]}).
		\end{align*}
		Since $\car(K)\neq 2$, we have the required equality.		
		
		The other equality is similar, using the commutative diagram
	
		\begin{center}
			\begin{tikzcd}
				{[[a,b],c]}\arrow[d,"{\tau_{[a,b],c}}"]\arrow[rr,"{[\tau_{a,b},e(c)]}"]& & {[[b,a],c]}\arrow[d,"{\tau_{[b,a],c}}"]\\
				{[c,[a,b]]}\arrow[rr,"{[e(c),\tau_{a,b}]}"]&  &{[c,[b,a]]}.
			\end{tikzcd}
		\end{center}
	\end{proof}
\end{prop}

\begin{defi}
	A \emph{braided internal functor between two braided categorical Lie $K$-algebras}, whose braidings are $\tau$ and $\tau'$,
is an internal functor $(F_1,F_0)$ such that $F_1(\tau_{a,b})=\tau'_{F_0(a),F_0(b)}$ for $a,b\in C_0$.
	
	We denote the category of braided categorical Lie $K$-algebras and braided internal functors between them as $\textbf{\textit{BICat}}(\textbf{\textit{LieAlg}}_K)$.
\end{defi}

The definition of braiding for crossed modules of Lie $K$-algebras was given in~\cite{TFM} trying to make a definition for which the Lie bracket
 was a braided for the crossed module $(M,M,[-,-],\Id_M)$ (the Lie bracket it is also known as commutator of the Lie $K$-algebra)
in parallelism to the fact that $(G,G,\Conj,\Id_G,[-,-])$ is a braided crossed module of groups. That definition
was also made to be a particular case of $2$-crossed modules of groups, whose definition can be seen in~\cite{M&P}.

Another definition can be seen in~\cite{Ulua}, but that definition does not verify the
mentioned requirements.

\begin{defi}
	Let $\mathcal{X}=(M,N,\cdot,\partial)$ be a crossed module of Lie $K$-algebras.
	
	A \emph{braiding} (or \emph{Peiffer lifting}) on the crossed module $\mathcal{X}$ is a $K$-bilinear map
	$\{-,-\}\colon N\times N\xrightarrow{ } M$ verifying:
	\begin{align}
	\partial\{n,n'\}&=[n,n'],\tag{BLie1}\label{BLie1}\\
	\{\partial m, \partial m' \}&=[m,m'],\tag{BLie2}\label{BLie2}\\
	\{\partial m, n \}&=-n\cdot m, \tag{BLie3}\label{BLie3}\\
	\{n,\partial m \}&=n\cdot m, \tag{BLie4}\label{BLie4}\\
	\{n,[n',n'']\}&=\{[n,n'],n''\}-\{[n,n''],n'\},\tag{BLie5}\label{BLie5}\\
	\{[n,n'],n''\}&=\{n,[n',n'']\}-\{n',[n,n'']\},\tag{BLie6}\label{BLie6}
	\end{align}
	for $m,m'\in M$, $n,n',n''\in N$.
	
	If $\{-,-\}$ is a braiding on $\mathcal{X}$ we will say that \emph{$(M,N,\cdot,\partial,\{-,-\})$ is a braided crossed module of Lie $K$-algebras}.
\end{defi}

\begin{example}
	It is immediate to check that $(M,M,[-,-],\Id_M,[-,-])$ is a braided crossed module of Lie $K$-algebras.
\end{example}

\begin{defi}
	An \emph{homomorphism of braided crossed modules of Lie $K$-algebras} $(M,N,\cdot,\partial,\{-,-\})\xrightarrow{(f_1,f_2)}(M',N',*,\partial',\{-,-\}')$
is an homomorphism of crossed modules of Lie $K$-algebras verifying $f_1(\{n,n'\})=\{f_2(n),f_2(n')\}'$ for $n,n'\in N$.
	
	We denote the category of  braided crossed modules of Lie $K$-algebras and its homomorphisms as $\textbf{\textit{BX}}(\textbf{\textit{LieAlg}}_K)$.
\end{defi}

Now, we will show the natural relation between definitions of braiding in the case of crossed modules of associative algebras and crossed modules of Lie algebras.

\begin{prop}
	If $\car(K)\neq 2$ and $(M,N,*,\partial,\{-,-\})$ is a braided crossed module of associative $K$-algebras,
 then $\{n,n'\}_{\mathfrak{L}}=\frac{\{n,n'\}-\{n',n\}}{2}$ is a braiding on the crossed module $(M^\mathcal{L},N^\mathcal{L},[-,-]_*,\partial)$.
	\begin{proof}
	We only need to show that the braiding axioms are verified. We will take $n,n',n''\in N, m,m'\in M$.
	
		The axioms \eqref{BLie1}, \eqref{BLie2}, \eqref{BLie3} and \eqref{BLie4} are proved using \eqref{BAs1}, \eqref{BAs2}, \eqref{BAs3} and \eqref{BAs4}.
		\begin{align*}
		\partial(\{n,n'\}_{\mathcal{L}})&=\partial(\frac{\{n,n'\}-\{n',n\}}{2})=\frac{[n,n']-[n',n]}{2}=[n,n'],\\	
		\{\partial m,\partial m'\}_{\mathcal{L}}&=\frac{\{\partial m,\partial m'\}-\{\partial m',\partial m\}}{2}=\frac{[m,m']-[m',m]}{2}=[m,m'],\\
		\{\partial m,n\}_{\mathcal{L}}&=\frac{\{\partial m,n\}-\{n,\partial m\}}{2}=\frac{-[n,m]_*-[n,m]_*}{2}=-[n,m]_*,\\
		\{n,\partial m\}_{\mathcal{L}}&=\frac{\{n,\partial m\}-\{\partial m,n\}}{2}=\frac{[n,m]_*+[n,m]_*}{2}=[n,m]_*.
		\end{align*}

		With the previous axioms proved, we have that the followings equalities hold.
		\begin{align*}
		\{n,[n',n'']\}_\mathcal{L}=\{n,\partial\{n',n''\}_\mathcal{L}\}_\mathcal{L}&=[n,\{n',n''\}_\mathcal{L}]_*\\&=-\{\partial\{n',n''\}_\mathcal{L},n\}_\mathcal{L}=-\{[n',n''],n\}_\mathcal{L}.
		\end{align*}
		Finally, we will prove that the braiding verifies the last axioms. For that we will abuse of language and we will denote $*=[-,-]_*$ (in the definition $*=(*_1,*_2)$).
		\begin{align*}
		&\{n,[n',n'']\}_\mathcal{L}=\{n,n'n''\}_\mathcal{L}-\{n,n''n'\}_\mathcal{L}\\
		&=\frac{1}{2}(\{n,n'n''\}-\{n'n'',n\}-\{n,n''n'\}+\{n''n',n\})\\
		&=\frac{1}{2}(n'*_1\{n,n''\}-\{n,n'\}*_2 n''-n'*_1\{n'',n\}-\{n',n\}*_2n''\\
		&-n''*_1\{n,n'\}-\{n,n''\}*_2n'+n''*_1\{n',n\}+\{n'',n\}*_2 n')\\
		&=\frac{1}{2}(n'*\{n,n''\}-n''*\{n,n'\}+n''*\{n',n\}-n'*\{n'',n\})\\
		&=n'*\{n,n''\}_\mathcal{L}-n''*\{n,n'\}_\mathcal{L}=-\{[n,n''],n'\}+\{[n,n'],n''\}.
		\end{align*}
		\begin{align*}
		& \{[n,n'],n''\}_\mathcal{L}=\{nn',n''\}_\mathcal{L}-\{n'n,n''\}_\mathcal{L}\\
		&=\frac{1}{2}(\{nn',n''\}-\{n'',nn'\}-\{n'n,n''\}+\{n'',n'n\})\\
		&=\frac{1}{2}(n*_1\{n',n''\}+\{n,n''\}*_2 n'-n*_1\{n'',n'\}-\{n'',n\}*_2n'\\
		&-n'*_1\{n,n''\}-\{n',n''\}*_2 n+n'*_1\{n'',n\}+\{n'',n'\}*_2 n)\\
		&=\frac{1}{2}(n*\{n',n''\}-n*\{n'',n'\}-n'*\{n,n''\}+n'*\{n''n\})\\
		&=n*\{n',n''\}_\mathcal{L}-n'*\{n,n''\}_\mathcal{L}=\{n,[n',n'']\}_\mathcal{L}-\{n',[n,n'']\}_\mathcal{L}.\qedhere
		\end{align*}
	\end{proof}
\end{prop}

\begin{remark}
	Note that the previous construction translates the example given for associative case to the one given in Lie case.
\end{remark}

\begin{remark}
	As for the case of groups we have, when $\car(K)\neq 2$, an equivalence between the categories
$\textbf{\textit{BICat}}(\textbf{\textit{LieAlg}}_K)$ and $\textbf{\textit{BX}}(\textbf{\textit{LieAlg}}_K)$. This can be seen in~\cite{TFM}.
	
	In addition, the relations with the associative case gives us two functors. This functors commute in a immediate way with the functors of the equivalence.
\end{remark}

\section{The non-abelian tensor product as example of braiding}\label{S:Lienonabtensor}

We will start with the non-abelian tensor product of groups. That product was introduced by Brown and Loday in~\cite{BrLod} through the following definition.

\begin{defi}
	Let $G$ and $H$ be two groups so that $G$ acts in $H$ with $\cdot$ an $H$ acts on $G$ with $*$, both by automorphisms.
	
	The \emph{non-abelian tensor product of $G$ with $H$} is denoted by $G\otimes H$ and is the group generated by the symbols $g\otimes h$, where $g\in G$, $h\in H$, and the relations
	\begin{align}
	gg'\otimes h&=(gg'g^{-1}\otimes g\cdot h)(g\otimes h),\tag{RTG1}\\
	g\otimes hh'&=(g\otimes h)(h*g\otimes hh'h^{-1}).\tag{RTG2}
	\end{align}
\end{defi}

The following proposition is mentioned for a general case in~\cite{BrLod}, using actions which are denominated compatible actions for make the tensor product.
 In this paper we will write a particular case, which is interesting to make an example.

\begin{prop}[\cite{BrLod}]
	Let $G$ be a group.
	
	Then $(G\otimes G,G,\cdot,\partial)$ is a crossed module of groups where $G\otimes G$ is the non-abelian tensor product of $G$ with itself using the conjugation action.
The action $\cdot\colon G\times (G\otimes G)\xrightarrow{}(G\otimes G)$ and the map $\partial \colon G\otimes G\xrightarrow{} G$ are defined on generators
 as $g\cdot(g_1\otimes g_2)=gg_1g^{-1}\otimes gg_2g^{-1}$ and $\partial(g_1\otimes g_2)=[g_1,g_2]$.
\end{prop}
This crossed module can be associated with a natural braiding, being this shown as an example in~\cite{Fuk}.
\begin{example}
	Let $G$ be a group.
	
	The map $\{-,-\}\colon G\times G\xrightarrow{} G\otimes G$ defined as $\{g_1,g_2\}=g_1\otimes g_2$ is a braiding over $(G\otimes G,G,\cdot,\partial)$, being this the one defined in the previous proposition.
	
	This is a immediate result, using the properties of the non-abelian tensor product of groups (see~\cite[Proposition 1.2.3]{McDer}) and the definition.
\end{example}

Once given the example in groups, we look for its analogue in Lie $K$-algebras. For this we need the concept of non-abelian tensor product of Lie $K$-algebras, introduced by Ellis in~\cite{Ellis}.
\begin{defi}
	Let $M$ and $N$ be two Lie $K$-algebras such that $M$ acts in $N$ by $\cdot$ and $N$ acts in $M$ with $*$.
	
	Its \emph{non-abelian tensor product} is denoted as  $M\otimes N$ and its defined as the Lie $K$-algebra generated by the symbols $m\otimes n$ with $m\in M$, $n\in N$ and the relations
	\begin{equation}\label{RTLie1}
	\lambda(m\otimes n)=\lambda m\otimes n=m\otimes \lambda n,\tag{RTLie1}
	\end{equation}
	\begin{align*}\label{RTLie2}
	(m+m')\otimes n=m\otimes n+m'\otimes n,\tag{RTLie2}\\
	m\otimes (n+n')=m\otimes n+m\otimes n',
	\end{align*}	\begin{align*}\label{RTLie3}
	[m,m']\otimes n=m\otimes (m'\cdot n)-m'\otimes (m\cdot n),\tag{RTLie3}\\
	m\otimes [n,n']=(n'*m)\otimes n-(n*m)\otimes n',
	\end{align*}
	\begin{equation}\label{RTLie4}
	[(m\otimes n),(m'\otimes n')]=-(n*m)\otimes(m'\cdot n'),\tag{RTLie4}
	\end{equation}
	where $m,m'\in M$, $n,n'\in N$.
\end{defi}

The following proposition, following the script of the case of groups, is proved more generally in the cite, but we restrict ourselves to the case that interests us.

\begin{prop}[\cite{Ellis}]
	Let $M$ be a Lie $K$-algebra.
	
	Then $(M\otimes M,M,\cdot,\partial)$ is a crossed modules of Lie $K$-algebras, where $M\otimes M$ is the non-abelian tensor product of $M$ with itself using the adjoin action.
 The action $\cdot\colon M\times (M\otimes M)\xrightarrow{}(M\otimes M)$ and the map $\partial \colon M\otimes M\xrightarrow{} M$ are defined
 on generators as $m\cdot(m_1\otimes m_2)=[m,m_1]\otimes m_2 +m_1\otimes[m,m_2]$ and $\partial(m_1\otimes m_2)=[m_1,m_2]$ where $[-,-]$ is the bracket of $M$.
\end{prop}

\begin{remark}
	We will rewrite, for clarity, the relations \eqref{RTLie3} and \eqref{RTLie4}  for the case of $M\otimes M$ with the
	 of $M$ in itself by adjoin action.
	\begin{align*}
	[m_1,m_2]\otimes m_3=m_1\otimes [m_2,m_3]-m_2\otimes [m_1,m_3],\tag{RTLie3}\\
	m_1\otimes [m_2,m_3]=[m_3,m_1]\otimes m_2-[m_2,m_1]\otimes m_3,
	\end{align*}
	\begin{equation}
	[(m_1\otimes m_2),(m_3\otimes m_4)]=[m_1,m_2]\otimes[m_3, m_4],\tag{RTLie4}
	\end{equation}
	where $m_1,m_2,m_3,m_4\in M$. For the last relation we use the anticommutativity.
\end{remark}

With this we will show an analogous example to the case of groups in the case of Lie $K$-algebras.

\begin{example}
	Let $M$ be a Lie $K$-algebra.
	
	The $K$-bilinear map $\{-,-\}\colon M\times M\xrightarrow{} M\otimes M$ defined with the expression $\{m_1,m_2\}=m_1\otimes m_2$ is a braiding
over the crossed module of Lie $K$-algebras $(M\otimes M,M,\cdot,\partial)$.
	
	We will check this. We will start with \eqref{BLie1}. If $m,m'\in M$, then
	\begin{align*}
	\partial\{m,m'\}=\partial(m\otimes m')=[m,m'].
	\end{align*}
	
	To check \eqref{BLie2}, by the $K$-linearity and $K$-bilinearity, we will work on generators, since the general case is only a sum of theirs.
	
	If $m_1\otimes m_2$ and $m_3\otimes m_4$ are generators of $M\otimes M$, then
	\begin{align*}
	\{\partial (m_1\otimes m_2), \partial (m_3\otimes m_4) \}=\{[m_1,m_2],[m_3,m_4]\}&=[m_1,m_2]\otimes[m_3,m_4]\\&=[(m_1\otimes m_2),(m_3\otimes m_4)],
	\end{align*}
	where the last equality is given by \eqref{RTLie4}.
	
	For the following properties we need a previous result.
	
	We will use \eqref{RTLie3} to prove $m_1\otimes [m_2,m_3]=-[m_2,m_3]\otimes m_1$.
	\begin{align*}
	&[m_1,m_2]\otimes m_3=m_1\otimes [m_2,m_3]-m_2\otimes [m_1,m_3]\\
	&=m_1\otimes [m_2,m_3]-[m_3,m_2]\otimes m_1+[m_1,m_2]\otimes m_3.
	\end{align*}
	Simplifying we have $0=m_1\otimes [m_2,m_3]+[m_2,m_3]\otimes m_1$, which is the wished equality.
	
	Now, we will show \eqref{BLie3}. Let $m\in M$ and $m_1\otimes m_2\in M\otimes M$.
	\begin{align*}
	&\{\partial (m_1\otimes m_2), m \}=\{[m_1,m_2],m\}=[m_1,m_2]\otimes m\\
	&=m_1\otimes [m_2,m]-m_2\otimes [m_1,m]=-m_1\otimes [m,m_2]+m_2\otimes [m,m_1]\\
	&=-m_1\otimes [m,m_2]-[m,m_1]\otimes m_2=-m\cdot (m_1\otimes m_2),
	\end{align*}
	where we use \eqref{RTLie3} together with the previous result.
	
	Now, we will verify \eqref{BLie4}. 	
	\begin{align*}
	\{m,\partial (m_1\otimes m_2) \}=m\otimes [m_1,m_2]&=-[m_1,m_2]\otimes m=- \{\partial(m_1\otimes m_2),m \}\\
	&=-(-m\cdot (m_1\otimes m_2))=m\cdot (m_1\otimes m_2),
	\end{align*}
	where we use \eqref{BLie3} and $m\otimes [m_1,m_2]=-[m_1,m_2]\otimes m$.
	
	Now, we will verify \eqref{BLie5} and \eqref{BLie6}. Let $m,m',m''\in M$.
	\begin{align*}
	&\{m,[m',m'']\}=m\otimes [m',m'']=[m'',m]\otimes m'-[m',m]\otimes m''=\\
	&=[m,m']\otimes m''-[m,m'']\otimes m'=\{[m,m'],m''\}-\{[m,m''],m'\},\\
	&\{[m,m'],m''\}=[m,m']\otimes m''=m\otimes [m',m'']-m'\otimes [m,m'']\\
	&=\{m,[m',m'']\}-\{m',[m,m'']\}.
	\end{align*}
	We use \eqref{RTLie3} in the second equality of both chains of equalities.
	
	Then, it is shown that $\{m,m'\}=m\otimes m'$ is a braiding.
\end{example}

\begin{remark}
	Note that the action given in the previous example is actually given as:
	\begin{align*}
	m\cdot (m_1\otimes m_2)=m\otimes [m_1,m_2].
	\end{align*}
\end{remark}

%

\end{document}